\documentclass[12pt]{amsart}

\usepackage{graphicx}
\usepackage{geometry}
\usepackage{latexsym}
\usepackage{amssymb}
\usepackage{amsmath}
\usepackage{cite}
\usepackage{bm}
\newtheorem{theorem}{Theorem}
\newtheorem{Remark}{Remark}
\newtheorem{lemma}{Lemma}

\newtheorem{proposition}{Proposition}
\newtheorem{definition}{Definition}
\newtheorem{question}{Question}

\numberwithin{equation}{section}

\begin{document}

\title{Mixed $L_p$ projection inequality }
\thanks{Supported by the Fundamental Research Funds for the Central Universities(2015ZCQ-LY-01, 2017ZY44)}

\author{Zhongwen Tang}
\address{College of Science\\ Beijing Forestry University\\ Beijing 100083, P.R.China}
\email{tzwlxsx@163.com}

\author{Lin Si$^*$}
\address{College of Science\\ Beijing Forestry University\\ Beijing 100083, P.R.China}
\email{silincd@163.com}
\thanks{$^*$ Corresponding author}

\begin{abstract}
In this paper, the mixed $L_{p}$-surface area measures are defined and the mixed $L_{p}$ Minkowski inequality is obtained consequently. Furthermore, the mixed $L_{p}$ projection inequality for mixed projection bodies is established.

\vskip 1ex

\end{abstract}

\maketitle

\section{Introduction}

In the $n$-dimensional Euclidean space $\mathbb{E}^{n}$, mixed volumes, which are the generalization of surface, volume, mean width and so on of a convex body(compact, convex subsets with nonempty interiors in
$\mathbb{E}^{n}$) are the core concept in the classical Brunn-Minkowski theory.
Surface area measures are one of important geometric measures of convex bodies in $\mathbb{E}^{n}$ for the integral formula of mixed volumes.
Usually, the surface area measure is  interpreted as the first variation of volume with respect to the Minkowski addition.

The mixed volume $V_1(K, L)$ of convex bodies $K$ and $L$, admits the following integral representation

$$V_{1}(K,L)=\frac{1}{n} \lim_{\varepsilon\rightarrow 0^+}\frac{V(K+\varepsilon L)-V(K)}{\varepsilon}=\frac{1}{n} \int_{S^{n-1}}h_LdS_K(u),$$
where $h_L$ is the support function of $L$, and $S_K(u)$ is the area measure of $K$.
The famous Minkowski's first inequality states that
\begin{equation}
V_{1}^{n}(K,L)\geq V^{n-1}(K)V(L),           \label{MFI}
\end{equation}
with equality if and only if $K$ and $L$ are homothetic.

In {\cite{firey62}}, Firey 
introduced the $L_p$ Minkowski combination (also known as the Minkowski-Firey combination) for each $p\geq 1$.
The $L_p$ Brunn-Minkowski theory, combined the $L_p$ Minkowski combination and the volume as a generalization of the classical Brunn-Minkowski theory has attracted increasing interest in recent years (see, e.g., \cite{lutwak93,lyz18,lh16,lyz00,lyz05,lutwak96,zx14,lhx18,lyz04,lv19,cg02}). One of the most important concept in the $L_p$ Brunn-Minkowski theory (see, {\cite{lutwak93}}) is the $L_p$-surface area measure which is a Borel measure (defined on the unit sphere $S^{n-1}$ in $\mathbb{E}^{n}$) for each convex body in $\mathbb{E}^{n}$ that contains the origin in its interior, and has a representation formula for the $p$-mixed volume due to Lutwak. The $L_{p}$-surface area measure and its associated Minkowski problem in the $L_{p}$ Brunn-Minkowski theory were introduced in {\cite{lutwak93}}.

Let $\mathcal{K}^{n}$ denote the class of convex bodies  in $\mathbb{E}^{n}$ and $\mathcal{K}_{o}^{n}$ denote the subset of $\mathcal{K}^{n}$ that contains the origin as interiors.

If $K,L\in \mathcal{K}_{o}^{n}$ and $p\geq 1$, then the $L_p$-mixed volume $V_{p}(K,L)$ is defined by
$$V_{p}(K,L)=\frac{p}{n}\lim_{\varepsilon\rightarrow 0^+}\frac{V(K+_p\varepsilon\cdot L)-V(K)}{\varepsilon}=\frac{1}{n}\int_{S^{n-1}}h_LdS_p(K,u),$$
where $S_{p}(K,\cdot)$ is the $L_{p}$-surface area measure of $K$. The $L_{p}$-mixed volume $V_{p}(K,L)$ is the case $i=0$ of the $L_{p}$-mixed quermassintegral introduced by Lutwak in {\cite{lutwak93}}.

The $L_{p}$ Minkowski inequality states that for $p\geq 1$,
\begin{equation}
V_{p}^{n}(K,L)  \geq V^{n-p}(K) V^{p}(L),  \label{LPMI}
\end{equation}
with equality if and only if $K$ and $L$ are dilates.
Obviously, the case $p=1$ of (\ref{LPMI}) is the Minkowski's first inequality (\ref{MFI}).

If $K,L\in \mathcal{K}_{o}^{n}$, for $p\geq 1$ and $0\leq i \leq n-1$, the $L_{p}$-mixed quermassintegrals  $W_{p,i}(K,L)$ is defined by
\begin{equation}
 W_{p,i}(K,L)=\frac{p}{n-i}\mathop{\lim}_{\varepsilon\rightarrow 0^+} \frac{W_{i}(K+_{p} \varepsilon \cdot L)-W_{i}(K)}{\varepsilon}. \label{lpmq}
\end{equation}
The Minkowski inequality for the $L_{p}$-mixed quermassintegrals states that for $p\geq 1$ and $0\leq i \leq n-1$
\begin{equation}
W_{p,i}^{n-i}(K,L)\geq W_{i}^{n-i-p}(K)W_{i}^{p}(L), \label{milpmq}
\end{equation}
with equality if and only if $K$ and $L$ are dilates.

The purpose of this paper is to construct mixed $L_p$-surface area measures and obtain the mixed $L_{p}$ projection inequalities.

\begin{definition}
For $K,L\in \mathcal{K}_{o}^{n}$, $Q\in \mathcal{K}^{n}$ and $0\leq t\leq n-1$, $p\geq 1$, the $L_{p,t}$ mixed volume is defined by
$$V_{p,t}(K,L,Q)
=\frac{ p}{t+1}\mathop{\lim}_{\varepsilon \rightarrow 0^{+}} \frac{V(K+_{p} \varepsilon \cdot L,t+1;Q,n-t-1)-V(K,t+1;Q,n-t-1) }{\varepsilon}. $$
\end{definition}

Our first main result is the following $L_{p}$-type Minkowski inequality for $L_{p,t}$ mixed volume.

\begin{theorem}
 If $K,L\in \mathcal {K}_{o}^{n}, Q\in \mathcal{K}^{n}$, then for $p\geq 1$ and $0\leq t\leq n-1$
\begin{equation}
 V_{p,t}^{n}(K,L,Q)\geq V^{{t+1-p}}(K) V^{{p}}(L) V^{{n-t-1}}(Q),   \label{MLPMI}
 \end{equation}
with equality if and only if $K$, $L$ are dilates and $Q$ is dilate, up to translation.
 \end{theorem}

An important inequality involving the polar projection bodies was obtained by Petty {\cite{petty72}} (see \cite{lutwak85} for an alternate proof). This inequality is now known as the Petty projection inequality, i.e.,
for $K\in \mathcal{K}^{n}$, then
$$V^{n-1}(K)V(\Pi^{*} K)\leq \omega_{n}^{n},$$
with equality if and only if $K$ is an ellipsoid. Here $\Pi K$ is the projection body of $K$ and $\Pi^{*} K$ is the polar body of $\Pi K$.

More generally, Lutwak (see, {\cite{lutwak85} or \cite{lutwak86}}) obtained the following mixed projection inequality, i.e.,
for $K_{1},\cdots,K_{n-1} \in \mathcal{K}^n$, then
$$V(K_{1})\cdots V(K_{n-1})V(\Pi^{*}(K_{1},\cdots,K_{n-1}))\leq \omega_{n}^{n},$$
with equality if and only if the $K_{i}$ are homothetic ellipsoids.

In {\cite{lyz00}}, Lutwak, Yang and Zhang introduced $L_{p}$-projection
body $\Pi_{p}K$ and established the $L_{p}$-Petty projection inequality, i.e.,
for $K \in \mathcal{K}_{o}^{n}$ and $1\leq p<\infty$,
$$V(K)^{(n-p)/p}V(\Pi_{p}^{*} K)\leq \omega_{n}^{n/p},$$
with equality if and only if $K$ is an ellipsoid centered at the origin.
An alternative proof of the $L_{p}$-Petty projection inequality was given by Campi and Grochi{\cite{cg02}}.

In Section 2, we shall introduce the mixed $L_{p}$-projection body $\Pi_{p,t}(K,Q)$, which is an extension of $L_{p}$-projection body. The case $t=n-1$ of the mixed $L_{p}$-projection body is the $L_{p}$-projection body $\Pi_{p}K$. The second main result of this paper is the mixed $L_{p}$-projection inequalities, which is an extension of $L_{p}$-Petty projection inequality.

\begin{theorem}
If $K\in \mathcal {K}_{o}^{n}$, $Q\in \mathcal{K}^{n}$, then for $p>1$ and $0< t < n-1$,
\begin{equation}
V^{(t+1-p)/p}(K) V^{(n-t-1)/p}(Q)V(\Pi_{p,t}^{*}(K,Q))\leq \omega_{n}^{n/p},\label{VPTI}
\end{equation}
with equality if and only if $K$ and $Q$ are dilate ellipsoids centered at the origin.
\end{theorem}

\section{Preliminary results}

For $x\in \mathbb{E}^{n}  \backslash  \{0\}$, the support function of $K\in \mathcal{K}^{n}$ is defined by
\begin{equation}
h_{K}(x)=\mathop{\max}\{x\cdot y : y\in K  \},   \label{SP}
\end{equation}
where $x\cdot y$ denotes the standard inner product of $x$ and $y$ in $\mathbb{E}^{n}$. For $x\in \mathbb{E}^{n}$ and $\phi \in GL(n) $, the support function of the image $\phi K=\{  \phi x : x \in K \}$ is given by
\begin{equation}
h_{\phi K}(x)=h_{K}(\phi^{t}x)   \label{HPHI}
\end{equation}
where $\phi^{t}$ denotes the transpose of $\phi$. For $K,L\in \mathcal{K}^{n}$ and $s,t \geq 0$, the Minkowski combination $sK+tL$ is defined by
$$sK+tL=\{ sx+ty: x\in K, y\in L \},$$
or equivalently,
$$h_{sK+tL}=sh_{K}+t h_{L}.$$

For $K, L \in \mathcal{K}_{o}^{n}$ and $s,t\geq 0$, the $L_p$ Minkowski combination, $s\cdot K +_{p} t\cdot L$, for each $p\geq 1$, is defined by
$$h_{s\cdot K+_{p} t\cdot L}^{p} = s h_{K}^{p} + t h_{L}^{p}.$$

The mixed volume of $K_1,\cdots,K_{n} \in \mathcal{K}^n$ is denoted by $V(K_1,\cdots,K_{n})$  and is uniquely determined by the requirement that it be symmetric in its arguments.
 If $K,L,Q\in \mathcal{K}^{n}$, then $V(K,s;L,t;Q,n-s-t)$ will be used to  denote the mixed volume
$$V(K,\cdots,K,L,\cdots,L,Q,\cdots,Q),$$
 in which $K$ appears $s$ times, $L$ appears $t$ times and $Q$ appears $n-s-t$ times.

For $K\in \mathcal{K}_{o}^{n}$, one has $K+_{p} \varepsilon \cdot K=(1+\varepsilon)^{\frac{1}{p}}K$ by the definition of the $L_p$ Minkowski combination. Hence, from the Definition 1, we can get that
 $$V_{p,t}(K,K,Q)=V(K,t+1;Q,n-t-1)$$
 and
 $$V_{p,t}(K,K,K)=V(K)$$
 for all $p\geq 1$, $0 \leq t \leq n-1$.

The mixed area measure $S(K_{1},\cdots,K_{n-1};\cdot)$ associated with the convex bodies $K_{1},\cdots,K_{n-1}$ in $\mathbb{E}^{n}$ is a unique (positive) Borel measure on $S^{n-1}$, with the property that for a convex body $K$ in $\mathbb{E}^{n}$, one has the integral representation
$$V(K_1,\cdots,K_{n-1},K)=\frac{1}{n} \int_{S^{n-1}} h_{K}(u) dS(K_1,\cdots,K_{n-1};u).$$
If $K,Q\in \mathcal{K}^{n}$, then $S(K,t;Q,n-t-1;\cdot)$ will be used to denote the mixed surface area measures
$$S(K,\cdots,K,Q\cdots,Q;\cdot),$$
in which $K$ appears $t$ times and $Q$ appears $n-t-1$ times. Notice that $S(K,\cdots,K;\cdot)$=$S_{K}(\cdot)$.

If $K\in \mathcal{K}^{n}$ and $u \in S^{n-1}$, let $K^{u}$ denote the orthogonal projection of $K$ onto the 1-codimensional space $u^{\perp}$ perpendicular to $u$. The $(n-1)$-dimensional volume of the orthogonal projection, $\upsilon(K^{u})$, also called the brightness of $K$ in the direction $u$. If $K_{1},\cdots,K_{n-1}\in \mathcal{K}^{n}$ and $u\in S^{n-1}$, then the mixed volume of $K_{1}^{u},\cdots,K_{n-1}^{u}$ in $u^{\perp}$ is
written as $\upsilon (K_{1}^{u},\cdots,K_{n-1}^{u})$ and called the mixed brightness of $K_{1},\cdots.K_{n-1}$ in the direction $u$.

The projection body $\Pi K$ of $K\in \mathcal{K}^{n}$ is an origin-symmetric convex body defined by
$$h_{\Pi K}(u)=\upsilon (K^{u})=\frac{1}{2}\int_{S^{n-1}}|u\cdot v|dS_{K}(v), ~~u\in S^{n-1}.$$
More generally, if $K_1,\cdots,K_{n-1} \in \mathcal{K}^{n}$, then the mixed projection body of $K_1,\cdots,K_{n-1}$, $\Pi(K_{1},\cdots,K_{n-1})$, is an origin-symmetric convex body defined by (see, e.g., {\cite{lutwak85,lutwak86}})
$$h_{\Pi(K_1,\cdots,K_{n-1})}(u)=\upsilon (K_{1}^{u},\cdots,K_{n-1}^{u})=\frac{1}{2}\int_{S^{n-1}} |u\cdot v|dS(K_1,\cdots,K_{n-1};v), ~~~u\in S^{n-1}.$$
Obviously, $\Pi(K,\cdots,K)=\Pi K$.

Let $\omega_{n}=\pi^{n/2}/\Gamma(1+\frac{n}{2})$  denote the volume of unit ball $B$ in $\mathbb{E}^{n}$ for $n\geq 0$.
If $K\in \mathcal{K}_{o}^{n}$, the $L_{p}$ projection body $\Pi_{p}K$ ($p\geq 1$) of $K$ is an origin-symmetric convex body defined by{\cite{lyz00}
\begin{equation}
h_{\Pi_{p}K}^{p}(u)= \frac{1}{n \omega_{n}c_{n-2,p}} \int _{S^{n-1}} |u\cdot v|^{p}  dS_{p}(K,v),    ~~~u\in S^{n-1},    \label{LPPJ}
\end{equation}
where $c_{n,p}=\omega_{n+p}/\omega_{2}\omega_{n}\omega_{p-1}$ such that $\Pi_{p}B=B$.

For $p\geq 1$ and $i=0,1,\cdots,n-1$, the $L_{p}$-mixed projection body of $K\in \mathcal{K}^{n}$ , $\Pi_{p,i}K$, is an origin-symmetric convex body introduced by Wang and Leng \cite{wl07}.

If $K\in \mathcal{K}_{o}^{n}$, $Q\in \mathcal{K}^{n}$, for $p\geq 1$ and $0 \leq t \leq n-1$, the mixed $L_{p}$ projection body is defined by
\begin{equation}
h_{\Pi_{p,t}{(K,Q)}}^{p}(u)=\frac{1}{n \omega_{n} c_{n-2,p}} \int_{S^{n-1}} |u\cdot v|^{p} dS_{p,t}(K,Q;v),        \label{LPTP}
\end{equation}
where $S_{p,t}(K,Q;\cdot)$ is the mixed $L_{p}$-surface area measure of $K,Q$ defined in Section 3. From (\ref{LPTP}), we have $\Pi_{p,t}(B,B)=B$.

If $K\in \mathcal{K}_{o}^{n}$, the polar body $K^{*}$ of $K$ is defined by
$$K^{*}=\{ x\in \mathbb{E}^{n}: x\cdot y \leq 1 ~~for~~~ all~~ y\in K  \}.$$
A compact set $K$ in $\mathbb{E}^{n}$ is star-shaped (about the origin) whose radial function is defined for $x\neq 0$ by
$$\rho_{K}(x)=max\{\lambda\geq 0: \lambda x \in K  \}.$$

Obviously, for $x\neq 0$ and $\phi \in SL(n)$, $\rho_{\phi K}(x)=\rho_{K}(\phi^{-1} x)$. For $K\in \mathcal{K}^{n}$, we can get that
\begin{equation}
\rho_{K^{*}}=1/h_{K}, ~~ and ~~ h_{K^{*}}=1/\rho_{K}.  \label{HRE}
\end{equation}

From (\ref{HRE}), for $K\in \mathcal{K}^{n}$ and $\phi \in SL(n)$, we have
\begin{equation}
(\phi K)^{*}=\phi^{-t} K,    \label{PK}
\end{equation}
where $\phi^{-t}$ denotes the inverse of the transpose of $\phi$.

If $\rho_{K}$ is positive and continuous, $K$ is called a star body (about the origin). Let $\mathcal{S}_{o}^{n}$ denote the set of star bodies (about the origin) in $\mathbb{E}^{n}$.
For $K,L\in \mathcal{S}_{o}^{n}$, $p\geq 1$ and $\varepsilon > 0$, the $L_{p}$-harmonic radial combination $K+_{-p} \varepsilon \cdot L \in \mathcal{S}_{o}^{n}$ is defined by (see {\cite{lutwak96}})
$$\rho(K+_{-p} \varepsilon \cdot L,\cdot)^{-p}=\rho_{K}(\cdot)^{-p} + \varepsilon \rho _{L}(\cdot)^{-p}.$$

If $K,L\in \mathcal{S}_{o}^{n}$, for $p\geq 1$, the $L_{p}$-dual mixed volume, $\tilde {V}_{-p}(K,L)$, of the $K$ and $L$ is defined by
$$ \tilde {V}_{-p}(K,L)=\frac{-p}{n}\mathop{\lim}_{\varepsilon\rightarrow o^{+}} \frac{V(K+_{-p} \varepsilon \cdot L )-V(K)}{\varepsilon}.$$
By the polar coordinate formula, one can get the following integral representation of the $L_{p}$-dual mixed volume $\tilde {V}_{-p}(K,L)$:
\begin{equation}
\tilde {V}_{-p}(K,L)=\frac{1}{n}\int_{S^{n-1}} \rho_{K}^{n+p}(v) \rho_{L}^{-p}(v) dS(v).   \label{VD}
\end{equation}

For $\phi \in SL(n)$, one can get that $\tilde {V}_{-p}(\phi K,\phi L)=\tilde {V}_{-p}(K,L)$ or equivalently
\begin{equation}
                 \tilde {V}_{-p}(\phi K,L)=\tilde {V}_{-p}(K,\phi^{-1}L).   \label{DVI}
\end{equation}
Obviously, for each star body $K$,
\begin{equation}
\tilde{V}_{-p}(K,K)=V(K).     \label{DI}
\end{equation}

The dual $L_{p}$ Minkowski mixed volume inequality states that
\begin{equation}
\tilde {V}_{-p}(K,L)\geq V(K)^{(n+p)/n} V(L)^{-p/n},   \label{VI}
\end{equation}
with equality if and only if $K$ and $L$ are dilates.
From the dual $L_{p}$ Minkowski mixed volume inequality (\ref{VI}) and identity (\ref{DI}), if $K_{1},K_{2}\in \mathcal{S}_{o}^{n}$, 
\begin{equation}
\tilde{V}_{-p}(L,K_{1})=\tilde{V}_{-p} (L,K_{2})    \label{V-P}
\end{equation}
for all star bodies $L$ which belong to some class that contains both $K_{1}$ and $K_{2}$, then $K_{1}=K_{2}$.

For general references on convex bodies and the $L_{p}$ Brunn-Minkowski theory, we refer to Gruber \cite{gruber} and Schneider \cite{schneider}.

\vskip 20pt

\section{Mixed $L_{p}$ surface area measures }

The Aleksandrov-Fenchel inequality (see, e.g., \cite{lutwak86}), in the form most suitable for our purposes, states that
\begin{equation}
V^{s}(K,s+t;Q,n-s-t)V^{t}(L,s+t;Q,n-s-t)\leq V^{s+t}(K,s;L,t;Q,n-s-t).      \label{AFI}
\end{equation}
 By repeated applications of the Aleksandrov-Fenchel inequalities and final application of the Minkowski inequality, Lutwak got the following result (see, also {\cite{schneider}, p.398}).
\begin{lemma}
If $K_{i} \in \mathcal{K}^{n}$, then
\begin{equation}
V^{n}(K_1,\cdots,K_{n})\geq V(K_1) \cdots V(K_{n}),\label{MLI}
\end{equation}
with equality if and only if the $K_{i}$ are homothetic.
\end{lemma}
For Lemma 1, in {\cite{adm99}}, Alesker, Dar and Milman gave an alternative proof in which ideas from mass transportation were used.

To establish the Theorem 1, the following result will be needed. The proof is similar to that of {\cite{lutwak93}}.

\vskip 10pt

\begin{theorem}
If $K,L\in \mathcal{K}_{o}^n$, $Q\in \mathcal{K}^{n}$, then for $p\geq 1$ and $0\leq t \leq n-1$, we have
\begin{align*}
&\mathop{\lim}_{\varepsilon \rightarrow 0^{+}} \frac{V(K+_{p} \varepsilon \cdot L,t+1;Q,n-t-1)-V(K,t+1;Q,n-t-1) }{\varepsilon}\\
&= \frac{t+1}{np} \int_{S^{n-1}} h_{L}^{p}(u) h_{K}(u)^{1-p} dS(K,t;Q,n-t-1;u).
\end{align*}
\end{theorem}

\vskip  10pt

\noindent{\bf Proof :}

Let
$$K_{\varepsilon}=K+_{p} \varepsilon \cdot L,$$
and define $g:[0,\infty) \rightarrow (0,\infty)$ by
$$g(\varepsilon)=V^{\frac{1}{t+1}}(K_{\varepsilon},t+1;Q,n-t-1).$$
Let
$$\mathop{\lim \inf}_{\varepsilon \rightarrow 0^+} \frac{V(K_{\varepsilon},t;K_{\varepsilon};Q,n-t-1)-V(K_{\varepsilon},t;K;Q,n-t-1)}{\varepsilon}=L_{\inf},$$
and
$$\mathop{\lim \sup}_{\varepsilon \rightarrow 0^+} \frac{V(K,t;K_{\varepsilon};Q,n-t-1)-V(K,t;K;Q,n-t-1)}{\varepsilon}=L_{\sup}.$$
Since $K_{\varepsilon}\supset K$, it follows from the monotonicity of the mixed volume that
$$L_{\inf}\geq 0 \  \    {\rm and}  \ \  ~ L_{\sup}\geq 0,$$
for all $\varepsilon$.
We shall make use of the fact that if $f_{0},f_{1},\cdots \in C(S^{n-1})$, with $\mathop{\lim}_{i\rightarrow \infty} f_{i} = f_{0}$, uniformly on $S^{n-1}$, and $ \nu_{0} ,\nu_{1},\cdots$ are finite measures on $S^{n-1}$ such that $\mathop {\lim}_{i\rightarrow \infty} \nu_{i}=\nu_{0},$ weakly on $S^{n-1}$, then
$$\mathop{\lim}_{i \rightarrow \infty} \int_{S^{n-1}} f_{i}(u) d\nu_{i}(u)=\int_{S^{n-1}} f_{0}(u) d \nu_{0}(u).$$
From the definition of the $L_{p}$ Minkowski combination, we have
$$\mathop{\lim}_{\varepsilon\rightarrow 0^+} \frac{h_{K_{\varepsilon} } -h_{K}}{\varepsilon}=\mathop{\lim}_{\varepsilon\rightarrow 0^+} \frac{(h_{K}^{p} + \varepsilon h_{L}^{p})^{\frac{1}{p}}-h_{K}}{\varepsilon}= \frac{1}{p} h_{L}^{p} h_{K}^{1-p},$$
uniformly on $S^{n-1}$.
By the weak continuity of mixed surface area measures $S(K_1,\cdots,K_{n-1},\cdot)$ (see, Schneider \cite{schneider}, p281) and the fact that $ \mathop{\lim}_{\varepsilon \rightarrow 0^+} K_{\varepsilon}=K$ in 
$\mathcal{K}^{n}$  (see,{\cite{firey62}}), we have
\begin{align*}
&\mathop{\lim }_{\varepsilon \rightarrow 0^+} \frac{V(K_{\varepsilon},t;K_{\varepsilon};Q,n-t-1)-V(K_{\varepsilon},t;K;Q,n-t-1)}{\varepsilon}\\
&=\mathop{\lim}_{\varepsilon \rightarrow 0^+} \frac{1}{n} \int_{S^{n-1}} \frac{h_{K_{\varepsilon}}-h_{K}}{\varepsilon}dS(K_{\varepsilon},t;Q,n-t-1;u)\\
&=\frac{1}{np} \int_{S^{n-1}} h_{L}^{p}(u) h_{K}^{1-p}(u) dS(K,t;Q,n-t-1;u).
\end{align*}
Similarly,
\begin{align*}
&\mathop{\lim }_{\varepsilon \rightarrow 0^+} \frac{V(K,t;K_{\varepsilon};Q,n-t-1)-V(K,t;K;Q,n-t-1)}{\varepsilon}\\
&=\mathop{\lim}_{\varepsilon \rightarrow 0^+} \frac{1}{n} \int_{S^{n-1}} \frac{h_{K_{\varepsilon}}-h_{K}}{\varepsilon}dS(K,t;Q,n-t-1;u)\\
&=\frac{1}{np} \int_{S^{n-1}} h_{L}^{p}(u) h_{K}^{1-p}(u) dS(K,t;Q,n-t-1;u).
\end{align*}
Then, we have
\begin{equation}
L_{\inf}=L_{\sup}=\frac{1}{np} \int_{S^{n-1}} h_{L}^{p}(u) h_{K}^{1-p}(u) dS(K,t;Q,n-t-1;u). \label{LLIS}
\end{equation}

From the definitions of $L_{\inf}$, $L_{\sup}$ and the Aleksandrov-Fenchel inequality (\ref{AFI}), we have that
$$\mathop{\lim \inf}_{\varepsilon\rightarrow 0^+}V^{\frac{t}{t+1}}(K_{\varepsilon},t+1;Q,n-t-1) \frac{V^{\frac{1}{t+1}}(K_{\varepsilon},t+1;Q,n-t-1)-V^{\frac{1}{t+1}}(K,t+1;Q,n-t-1)}{\varepsilon}\geq L_{\inf},$$
$$\mathop{\lim \sup}_{\varepsilon\rightarrow 0^+}V^{\frac{t}{t+1}}(K,t+1;Q,n-t-1) \frac{V^{\frac{1}{t+1}}(K_{\varepsilon},t+1;Q,n-t-1)-V^{\frac{1}{t+1}}(K,t+1;Q,n-t-1)}{\varepsilon}\leq L_{\sup}.$$

By the continuity of $V(K_{1},\cdots,K_{n})$ for $K_i \in \mathcal{K}^{n}$, $i=1,\cdots,n$, $g$ is continuous at $0$ (see, Schneider \cite{schneider}, P280). And notice that
$\mathop {\lim}_{\varepsilon \rightarrow 0^+} K_{\varepsilon}=K$
in $\mathcal{K}_{o}^{n}.$ Thus, the above inequalities can be rewritten as
$$V^{\frac{t}{t+1}}(K,t+1;Q,n-t-1)\mathop{\lim \inf}_{\varepsilon\rightarrow 0^+} \frac{V^{\frac{1}{t+1}}(K_{\varepsilon},t+1;Q,n-t-1)-V^{\frac{1}{t+1}}(K,t+1;Q,n-t-1)}{\varepsilon}\geq L_{\inf},$$
and
$$V^{\frac{t}{t+1}}(K,t+1;Q,n-t-1)\mathop{\lim \sup}_{\varepsilon\rightarrow 0^+} \frac{V^{\frac{1}{t+1}}(K_{\varepsilon},t+1;Q,n-t-1)-V^{\frac{1}{t+1}}(K,t+1;Q,n-t-1)}{\varepsilon}\leq L_{\sup}.$$
Since $L_{\inf}=L_{\sup}$, then the last two inequalities will imply that $g$ is differentiable at $0$, and $g(0)^{t}g^{'}(0) =L_{\inf}=L_{\sup}.$ The differentiability of $g$ at $0$ implies that $g^{t+1}$ is differentiable at $0$, then we have
$$\mathop{\lim}_{\varepsilon\rightarrow 0^+} \frac{g^{t+1}(\varepsilon)-g^{t+1}(0)}{\varepsilon}=(t+1)g^{t}(0)\mathop{\lim}_{\varepsilon\rightarrow 0^+} \frac{g(\varepsilon)-g(0)}{\varepsilon},$$
or
$$\frac{1}{t+1}\mathop{\lim}_{\varepsilon\rightarrow 0^+} \frac{V(K_{\varepsilon},t+1;Q,n-t-1)-V(K,t+1;Q,n-t-1)}{\varepsilon}=L_{\inf}=L_{\sup}.$$
\qed

Suppose $K\in \mathcal{K}_{o}^{n}$ and $Q\in \mathcal{K}^{n}$, the mixed $L_{p}$ surface area measure of $K,Q$ is defined by
\begin{equation}
dS_{p,t}(K,Q;\cdot)=h_{K}^{1-p}(\cdot) d S(K,t;Q,n-t-1;\cdot).  \label{MSAM}
\end{equation}
The $L_{p,t}$ mixed volume of $K,L\in \mathcal{K}_{o}^{n},Q\in \mathcal {K}^{n}$ can also be defined as
\begin{equation}
V_{p,t}(K,L,Q)=\frac{1}{n} \int_{S^{n-1}} h_{L}^{p}(u)d S _{p,t} (K,Q;u).  \label{LPT3}
\end{equation}
When $t=n-1$ in (\ref{LPT3}), the $L_{p,t}$ mixed volume induces to the $L_{p}$ mixed volume
$$V_{p,n-1}(K,L,Q)=V_{p}(K,L)=\frac{1}{n} \int_{S^{n-1}}  h_{L}^{p}(u)  dS_{p}(K,u).$$
When $p=1$ in (\ref{LPT3}), the $L_{p,t}$ mixed volume induces to the mixed volume of $K,L,Q\in \mathcal{K}^{n}$
$$V_{1,t}(K,L,Q)=V(K,t;L ;Q,n-t-1).$$
When $Q=B$ in (\ref{LPT3}), the $L_{p,t}$ mixed volume induces to the $L_{p}$ mixed quermassintegrals $W_{p,i}(K,L)$, which has the following integral representation,
$$V_{p,t} (K,L,B)= W_{p,i}(K,L)=\frac{1}{n} \int_{S^{n-1}} h _{L}^{p}(u) dS_{p,i}(K,u),$$
where $i=n-t-1$.

\vskip 10pt

\noindent{\bf Proof of Theorem 1 :}

Case I: $t> 0$

By the H\"{o}lder inequality and (\ref{MLI}), we have
\begin{align*}
V_{p,t}(K,L,Q)&=\frac{1}{n} \int_{S^{n-1}} h_{L}^{p}(u) h_{K}^{1-p}(u) dS(K,t;Q,n-t-1;u)\\
&\geq V^{p}(K,t;L;Q,n-t-1) V^{1-p}(K,t+1;Q,n-t-1)\\
&\geq \Big( V^{\frac{t}{n}}(K) V^{\frac{1}{n}}(L)V^{\frac{n-t-1}{n}} (Q) \Big)^p  \Big(  V^{\frac{t+1}{n}}(K) V^{\frac{n-t-1}{n}}(Q)  \Big)^{1-p}\\
&=V^{\frac{t+1-p}{n}}(K) V^{\frac{p}{n}}(L) V^{\frac{n-t-1}{n}}(Q).
\end{align*}

To obtain the equality conditions, notice that there is equality in H\"{o}lder's inequality precisely when
\begin{equation}
V(K,t;L;Q,n-t-1)h_{K}=V(K, t+1;Q,n-t-1)h_{L},    \label{VH}
\end{equation}
almost everywhere, with respect to the measure $S(K, t+1;Q,n-t-1;\cdot)$, on $S^{n-1}$. Equality in inequality (\ref{MLI}) holds precisely when there exists an $x\in \mathbb{E}^{n}$ such that
\begin{equation}
V(K,t;L;Q,n-t-1)h_{K}(u)= x\cdot u + V(K, t+1;Q,n-t-1)h_{L}(u),  \label{VHx}
\end{equation}
for all $u\in S^{n-1}$ and there exists a $y\in \mathbb{E}^{n}$ such that
$$\lambda h_{Q}(u)=y\cdot u + h_{L}(u) ,\ \ \lambda>0$$
for all $u\in S^{n-1}$.
Since the support of the measure $S(K, t;Q,n-t-1;\cdot)$ cannot be contained in the great sphere of $S^{n-1}$ orthogonal to $x$. Hence, from ({\ref{VH}}) and ({\ref{VHx}}), we have $x=0$ and
$$V(K,t;L;Q,n-t-1)h_{K}=V(K, t+1;Q,n-t-1)h_{L} $$
everywhere.

Case II: $t=0$

Notice that $V(K,0;L;Q,n-1)=V_{1}(Q,L)$ and $V(K,1;Q,n-1)=V_{1}(Q,K)$. Thus, we have
\begin{align*}
V_{p,0}(K,L,Q)&=\frac{1}{n} \int_{S^{n-1}} h_{L}^{p}(u)h_{K}^{1-p}(u)dS_{Q}(u)\\
&\geq V_{1}^{p}(Q,L) V_{1}^{1-p}(Q,K)\\
&\geq \Big( V^{\frac{n-1}{n}}(Q) V^{\frac{1}{n} } (L) \Big)^{p} \Big(  V^{\frac{n-1}{n}} (Q) V^{\frac{1}{n}} (K)   \Big)^{1-p}\\
&=V^{\frac{1-p}{n}} (K) V^{\frac{p}{n}}(L) V^{\frac{n-1}{n}}(Q).
\end{align*}
The proof of the equality conditions is the same as Case I. \qed

The classical Minkowski's problem consists in retrieving $K$ from its surface area
measure, and it is well-known that it admits a unique solution up to translations.
Corresponding to the classical Minkowski's problem, for the mixed $L_{p}$-surface area measure $S_{p,t}(K,Q;\cdot)$, $K\in \mathcal{K}_{o}^{n}$, $Q\in \mathcal{K}^{n}$, we have the following conjecture.

\begin{question}
If $\mu$ is an positive Borel measure on $S^{n-1}$, which is not concentrated on a great sphere of $S^{n-1}$, and $0\leq t\leq n-1$, then
there exist unique $K\in \mathcal{K}_{o}^{n}$ and $Q\in \mathcal{K}^{n}$ such that
$$S_{p,t}(K,Q;\cdot)=\mu?$$
\end{question}

\vskip 20pt

\section{Mixed $L_p$ projection inequalities}

For $p\geq 1$ and $K\in \mathcal{S}_{o}^{n}$, the $L_{p}$-centroid body $\Gamma_{p}K$ is an origin-symmetric body whose support function is given by (see, e.g., {\cite{lyz00}})
$$h_{\Gamma_{p}K}^{p}(x)= \frac{1}{c_{n,p}V(K)} \int_{K} |x\cdot y|^{p} dy=\frac{1}{(n+p)c_{n,p}V(K)} \int_{S^{n-1}} |x\cdot v|^{p} \rho_{K}^{n+p} (v)dS(v),$$
for all $x\in \mathbb{E}^{n}$. The normalization is chosen so that, for the unit ball $B$ in $\mathbb{E}^{n}$, we have $\Gamma_{p}B=B$.
For the polar of $\Gamma_{p}K$, we will write $\Gamma_{p}^{*}K$.

A consequence of the definition of the $L_{p}$-centroid body $\Gamma_{p}K$, is that for $\phi \in SL(n)$,
\begin{equation}
\Gamma_{p}\phi K=\phi \Gamma_{p}K.  \label{GP}
\end{equation}
Obviously, from the definition of $\Gamma_{p}K$, for $\lambda >0$, we also have
\begin{equation}
\Gamma_{p} \lambda K = \lambda \Gamma_{p}K. \label{GPLK}
\end{equation}

For the proof of the Theorem 2, the following lemmas are needed.

\begin{lemma}
If $K\in \mathcal{K}_{o}^{n}$, $L\in \mathcal{S}_{o}^{n}$, $Q \in  \mathcal{K}^{n} $, for $p>1$ and $1<t<n-1$, then
\begin{equation}
V_{p,t}(K,\Gamma_{p}L, Q )=\frac{\omega_{n}}{V(L)} \tilde{V}_{-p}(L,\Pi_{p,t}^{*}(K,Q)).    \label{LPDE}
\end{equation}
\end{lemma}
\noindent{\bf Proof :}

From the integral representation (\ref{LPT3}), the definition of the $L_{p}$-centroid body of $L$, Fubini's theorem, (\ref{HRE}) and (\ref{VD}),
we have
\begin{align*}
V_{p,t}(K,\Gamma_{p}L,Q)&=\frac{1}{n}\int_{S^{n-1}} h_{\Gamma_{p}L}^{p}(u) dS_{p,t}(K,Q,u)\\
&=\frac{1}{n}\int_{S^{n-1}} \frac{1}{c_{n,p} V(L)} \int_{L} |u\cdot x|^{p} dx  dS_{p,t}(K,Q,u)\\
&=\frac{1}{n c_{n,p} V(L) (n+p)}\int_{S^{n-1}} \int_{S^{n-1}} |u\cdot v|^{p} \rho _{L}^{n+p} (v) dS(v) dS_{p,t}(K,Q,u)\\
&=\frac{1}{n c_{n,p} V(L) (n+p)}\int_{S^{n-1}} \int_{S^{n-1}} |u\cdot v|^{p} dS_{p,t}(K,Q,u)\rho _{L}^{n+p} (v) dS(v)\\
&=\frac{\omega_{n}}{V(L)} \int_{S^{n-1}} \rho_{L}^{n+p} (v) h_{\Pi _{p,t}(K,Q)}^{p} (v) d S(v)\\
&=\frac{\omega_{n}}{V(L)} \int_{S^{n-1}} \rho_{L}^{n+p} (v) \rho_{\Pi _{p,t}^{*}(K,Q)}^{-p} (v) d S(v)\\
&=\frac{\omega_{n}}{V(L)} \tilde{V}_{-p} (L,\Pi_{p,t}^{*}(K,Q) ).
\end{align*}
\qed

In {\cite{lyz00}}, Lutwak, Yang and Zhang obtained the following result for the volume of the $L_{p}$-centroid body.
\begin{lemma}
If $K\in \mathcal{S}_{o}^{n}$, for $p> 1$, then
\begin{equation}
V(\Gamma_{p}K)\geq V(K),  \label{CBI}
\end{equation}
with equality if and only if $K$ is an ellipsoid centered at the origin.
\end{lemma}

The following Lemma 4 can be obtained by using the similar methods in {\cite{lyz05}}. A proof is given in the Appendix.

\begin{lemma}

If $K,L\in \mathcal{K}_{o}^{n}$, $Q\in \mathcal{K}^{n}$, for $p > 1$, $0<t<n-1$ and $\phi\in SL(n)$, then
\begin{equation}
V_{p,t}(\phi K,L,\phi Q)=V_{p,t}(K,\phi^{-1}L,Q).  \label{PTPI}
\end{equation}
\end{lemma}

A property of the operator $\Pi_{p,t}$ is as follows.
\begin{lemma}
If $K\in \mathcal{K}_{o}^{n}$, $Q\in \mathcal{K}^{n}$, for $p> 1$, $0<t<n-1$ and $\phi \in SL(n)$, then
\begin{equation}
\Pi_{p,t}(\phi K,\phi Q)= \phi ^{-t} \Pi _{p,t}(K,Q),  \label{LPPK}
\end{equation}
where $\phi^{-t}$ denotes the inverse transpose of $\phi$.
\end{lemma}
\noindent{\bf Proof:}

From (\ref{LPDE}), (\ref{PTPI}), (\ref{GP}) and (\ref{DVI}), for each $L\in \mathcal{S}_{o}^{n}$, we have
\begin{align*}
\frac{\omega_{n} V_{-p}(L,\Pi_{p,t}^{*}(\phi K,\phi Q))}{V(L)}&=V_{p,t}(\phi K,\Gamma_{p}L,\phi Q)\\
&=V_{p,t}(K,\phi^{-1}\Gamma_{p}L,Q)\\
&=V_{p,t}(K,\Gamma_{p} \phi^{-1}L,Q)\\
&=\frac{\omega_{n}V_{-p} (\phi^{-1}L,\Phi_{p,t}^{*}(K,Q))}{V(L)}\\
&=\frac{\omega_{n}V_{-p}(L,\phi \Pi_{p,t}^{*}(K,Q))}{V(L)}.
\end{align*}
According to ({\ref{V-P}}), for all $L\in \mathcal{S}_{o}^{n}$,
$$\frac{\omega_{n} V_{-p}(L,\Pi_{p,t}^{*}(\phi K,\phi Q))}{V(L)}=\frac{\omega_{n}V_{-p}(L,\phi \Pi_{p,t}^{*}(K,Q))}{V(L)}$$
implies that
$$\Pi_{p,t}^{*}(\phi K, \phi Q)=\phi \Pi_{p,t}^{*}(K,Q).$$
And from ({\ref{PK}}), we have the desired equality.  \qed

If $K\in \mathcal{K}_{o}^{n}$ and $Q\in \mathcal{K}$,  then for $p > 1$, $0< t< n-1$ and $\lambda_{1},\lambda_{2}>0$, we have
\begin{equation}
\Pi_{p,t}(\lambda_{1}K,\lambda_{2}Q)=\lambda_{1}^{(t+1-p)/p}\lambda_{2}^{(n-t-1)/p}\Pi_{p,t}(K,Q). \label{LBDA}
\end{equation}

\vskip 10pt

\noindent{\bf Proof of Theorem 2:}

In (\ref{LPDE}), let $L=\Pi_{p,t}^{*}(K,Q)$ and notice that $V_{-p}(K,K)=K$, we can get$$\omega_{n}=V_{p,t}(K,\Gamma_{p}\Pi_{p,t}^{*}(K,Q),Q).$$
Combining with (\ref{MLPMI}) and (\ref{CBI}), we have
\begin{align*}
\omega_{n}&\geq V^{{t+1-p}/{n}}(K) V^{{p}/{n}} (\Gamma_{p} \Pi _{p,t}^{*}(K,Q))V^{{n-t-1}/{n}}(Q)  \\
&\geq V^{{t+1-p}/{n}}(K)V^{{n-t-1}/{n}}(Q)V^{{p}/{n}} (\Pi_{p,t}^{*}(K,Q)).
\end{align*}
Then the inequality (\ref{VPTI}) is obtained.

From the equality conditions of (\ref{MLPMI}), the equality holds in the first inequality of the above inequalities if and only if $K$ and $\Gamma_{p} \Pi_{p,t}^{*}(K,Q)$ are dilates and $Q$ is up to translation and dilate. And from the equality condition of (\ref{CBI}), the equality holds in the second inequality of the above inequalities if and only if $\Pi_{p,t}^{*}(K,Q)$ is an ellipsoid centered at the origin.
 It is easily seen that if $\Pi_{p,t}^{*}(K,Q)$ is an ellipsoid centered at the origin, then $\Pi_{p,t}(K,Q)$ is an ellipsoid. Hence, for $\lambda>0$ and $\phi \in SL(n)$,
  \begin{equation}
  \Pi_{p,t}(K,Q)=(\lambda \phi)^{-1} B.    \label{PIB}
  \end{equation}
Let $\lambda=\lambda_{1}^{(t+1-p)/p} \lambda_{2}^{(n-t-1)/p}$, from ({\ref{LPPK}}), ({\ref{LBDA}}) and (\ref{PIB}), we have
$$\Pi_{p,t}(\lambda_{1} \phi^{-t}K,\lambda_{2} \phi^{-t} Q)=B.$$
Therefore, we have that $K=\lambda_{1}^{-1}\phi^{t}B, Q=\lambda_{2}^{-1}\phi^{t}B$. Then, from (\ref{PIB}), (\ref{PK}), (\ref{GP}), (\ref{GPLK}) and the fact that $\Gamma_{p}B=B$, we can obtain
$\Gamma_{p} \Pi_{p,t}^{*}(K,Q)=\lambda \phi^{t}B$.

Thus, we conclude that the equality holds in inequality (\ref{VPTI}) if and only if $K$ and $Q$ are dilate ellipsoids centered at the origin.   \qed

\vskip  20pt
\section*{Appendix}
The proof of the Lemma 4 is similar to that of Corollary 1.3 in {\cite{lyz05}}. For $x\in \mathbb{E}^{n}$ and $x\neq 0$, let $\langle x \rangle=x/|x|$.

\begin{definition}
Given a measure $\mu(u)$ on $S^{n-1}$, for $p > 0$, $\phi \in GL(n)$ and $f\in C(S^{n-1})$, define the measure $\mu^{(p)}(\phi u)$ on $S^{n-1}$ by
$$\int_{S^{n-1}} f(u) d\mu^{(p)} (\phi u) = \int_{S^{n-1}} |\phi^{-1}u|^{p} f(\langle\phi^{-1}u\rangle)d\mu(u).$$
\end{definition}

If $K_{1},\cdots,K_{n}\in \mathcal{K}^{n}$ and $\phi \in SL(n)$, then (see, {\cite{schneider}}, p282)
\begin{equation}
\int_{S^{n-1}} h_{\phi K_{n}}(u) dS(\phi K_{1},\cdots, \phi K_{n-1},u)=\int_{S^{n-1}} h_{K_{n}}(u) dS( K_{1},\cdots, K_{n-1},u).\label{PHI}
\end{equation}

For each convex body $L$, it follows from Definition 2, the homogeneity of $h_{L}$, (\ref{HPHI}) and (\ref{PHI}) that
\begin{align*}
\int_{S^{n-1}} h_{L}(u)dS^{(1)}(K,t;Q,n-t-1;\phi^{t} u) &=\int_{S^{n-1}} |\phi^{-t}u|h_{L}(\langle \phi^{-t}u  \rangle) dS(K,t;Q,n-t-1;u)\\
&=\int_{S^{n-1}} h_{L} (\phi^{-t}u)dS(K,t;Q,n-t-1;u)\\
&=\int_{S^{n-1}} h_{L} (u)dS(\phi K,t;\phi Q,n-t-1;u).
\end{align*}
From the fact that if two Borel measures on $S^{n-1}$ are equal when integrated against support functions of convex bodies then the measures are identical (see,{\cite{lyz05}}).
Thus for $K,Q \in \mathcal{K}^{n}$ and $\phi \in SL(n)$, we have

\begin{equation}
dS(\phi K,t;\phi Q,n-t-1;u)=dS^{(1)} (K,t;Q,n-t-1;\phi^{t}u).  \label{MSA}
\end{equation}

\begin{proposition}
If $K\in \mathcal {K}_{o}^{n}, Q\in \mathcal{K}^{n}$ and $p>0$, then for $\phi \in SL(n)$,
\begin{equation}
dS_{p,t}(\phi K,\phi Q,u)=dS_{p,t}^{(p)}(K,Q,\phi^{t}u)   \label{DPT}
\end{equation}
\end{proposition}
\noindent{\bf Proof:}

For $f\in C(S^{n-1})$, from (\ref{MSAM}), (\ref{HPHI}), (\ref{MSA}), Definition 2, the homogeneity of $h_{K}$, we have
\begin{align*}
\int_{S^{n-1}} f(u) dS_{p,t}(\phi K,\phi Q,u)&=\int_{S^{n-1}} f(u) h_{K}^{1-p}(\phi ^{t} u) dS(\phi K,t;\phi Q,n-t-1;u)\\
&=\int_{S^{n-1}} f(u) h_{K}^{1-p}(\phi^{t} u) dS^{(1)}(K,t;Q,n-t-1;\phi^{t}u)\\
&=\int_{S^{n-1}} |\phi^{-t}u| f(\langle \phi^{-t} u \rangle ) h_{K}^{1-p}(\phi^{t}  \langle \phi^{-t} u  \rangle) dS(K,t;Q,n-t-1;u)\\
&=\int_{S^{n-1}} |\phi^{-t}u|^{p} f(\langle \phi^{-t} u \rangle ) h_{K}^{1-p}(u) dS(K,t;Q,n-t-1;u)\\
&=\int_{S^{n-1}} |\phi^{-t}u|^{p} f(\langle \phi^{-t} u \rangle ) dS_{p,t}(K,Q,u)\\
&=\int_{S^{n-1}}  f(u) dS_{p,t}^{(p)}(K,Q,\phi^{t}u).
\end{align*}
\qed

\vskip 10pt
\noindent{\bf Proof of Lemma 4:}

From (\ref{LPT3}), Proposition 1, Definition 2, the homogeneity of the support function, and (\ref{HPHI}), we have
\begin{align*}
V_{p,t} (\phi K, L, \phi Q)&=\frac{1}{n} \int_{S^{n-1}} h_{L}^{p}(u) dS_{p,t}(\phi K,\phi Q,u)\\
&=\frac{1}{n} \int_{S^{n-1}} h_{L}^{p} (u) dS_{p,t}^{(p)} (K,Q,\phi^{t}u)\\
&=\frac{1}{n} \int_{S^{n-1}}  |\phi^{-1}u|^{p}  h_{L}^{p} ( \langle  \phi^{-t} u   \rangle  ) dS_{p,t} (K,Q,u)\\
&=\frac{1}{n} \int_{S^{n-1}}   h_{L}^{p} (  \phi^{-t} u   ) dS_{p,t} (K,Q,u)\\
&=V_{p,t}(K,\phi^{-1}L,Q).
\end{align*}
\qed

\vskip  30pt

\def\cprime{$'$}
\providecommand{\bysame}{\leavevmode\hbox to3em{\hrulefill}\thinspace}
\raggedbottom

\end{document}